\documentclass[a4paper, 10pt]{article}
\usepackage{amssymb, amsmath, amsfonts}
\usepackage[english]{babel}
\usepackage{amsthm,amscd,amsfonts,latexsym,color,epsfig}

\begin{document}

\begin{center}
{\bf NEW ESTIMATE FOR THE MULTINOMIAL MITTAG-LEFFLER FUNCTION} 
\end{center}
\begin{center}

\centerline {{\bf Murat Mamchuev}
\footnote{
Institute of Applied Mathematics and Automation of KBSC of RAS, 
Shortanova str. 89-A, Nalchik, 360000, Kabardino-Balkar Republic, Russia, E-mail: mamchuev@rambler.ru}}
\end{center}

\bigskip

\begin{abstract}
In this paper, a new estimate is obtained for the multinomial Mittag-Leffler function.
This function was introduced by Yuri Luchko and Rudolfo Gorenflo as the fundamental solution 
of the ordinary differential equation of fractional discrete distributed order.
\end{abstract}

\vspace{5mm}

\vspace{5mm}

The multinomial Mittag-Leffler function is defined as [1]
$$E_{(\mu_1,\mu_2,...,\mu_n),\gamma}(z_1,z_2,...,z_n)$$
$$=\sum\limits_{k=0}^{\infty}
\sum\limits_{l_1+...+l_n=k}(k;l_1,...,l_n)\frac{z_1^{l_1}\cdot...\cdot z_n^{l_1}}
{\Gamma(\gamma+\sum_{i=1}^{n}\mu_i l_i)},$$
here $(k;l_1,...,l_n)$ denotes the multinomial coefficient
$$(k;l_1,...,l_n)=\frac{k!}{l_1 !\cdot...\cdot l_n !},
\quad k=\sum\limits_{i=1}^{n}l_i,
$$
where $l_i (i=1,...,n)$ are non-negative integers.

We also need the following definition of the Mittag-Leffler-type function [2]
$$E_{\mu, \nu}(z)=\sum\limits_{k=0}^{\infty}\frac{z^k}{\Gamma(\mu k+\nu)}.$$

{\bf Lemma 1.} {\it Let $0<\mu_1<\mu_2<...<\mu_n,$ $\gamma>0,$ then the estimate 
$$\left |E_{(\mu_1,\mu_2,...,\mu_n),\gamma}(z_1,z_2,...,z_n)\right |\leq C E_{\mu_1,\gamma}(|z_1|+|z_2|+...+|z_n|)$$   
holds, here 
$C=1+\frac{\Gamma(\gamma+\mu_1 n_0)}{\Gamma(\gamma)},$ 
$n_0\in {\mathbb N}$ is a number
such that} 
$$\mu_1 n_0<x_0-\gamma<\mu_1 (n_0+1), \quad x_0=\min\limits_{x>0}\Gamma(x).$$

{\bf Proof.}
The following relation 
\begin{equation}\label{eq1}
\sum\limits_{l_1+...+l_n=k}(k;l_1,...,l_n)\prod\limits_{i=1}^{n}z_i^{l_i}=
(z_1+...+z_n)^k,
\end{equation}
can be proved by the mathematical induction method.

Note that
$$\gamma+\sum\limits_{i=1}^{n}\mu_i l_i=
\gamma+\mu_j \sum\limits_{i=1}^{n}l_i+\sum\limits_{i=1}^{n}(\mu_i-\mu_j)l_i=
\gamma+\mu_j k+\sum\limits_{i=1}^{n}(\mu_i-\mu_j)l_i, \quad 0\leq j\leq n.$$
Consequently,
$$\gamma+\mu_1 k <\gamma+\sum\limits_{i=1}^{n}\mu_i l_i< \gamma+\mu_n k.$$
Sinse the function $\frac{1}{\Gamma(x)}$ has only one maximum on the interval $(0,\infty)$
which reached at the point $x_0=1,4616321...,$ and decreases on the interval $(x_0,\infty),$ then
there exists the number $n_0\in {\mathbb N}$ such that the inequality
\begin{equation}\label{eq2}
\frac{1}{\Gamma(\gamma+\sum_{i=1}^{n}\mu_i l_i)}<
\frac{1}{\Gamma(\gamma+\mu_1 k)}
\end{equation}
holds for every numbers $k>n_0.$
For $0\leq k\leq n_0$ we have the inequality
\begin{equation}\label{eq3}
\frac{1}{\Gamma(\gamma+\sum_{i=1}^{n}\mu_i l_i)}<
C_0\frac{1}{\Gamma(\gamma+\mu_1 k)},
\end{equation}
where
$$C_0=\frac{\max\limits_{0\leq k\leq n_0}\Gamma(\gamma+\mu_1 k)}
{\min\limits_{0\leq k\leq n_0}\Gamma(\gamma+\sum_{i=1}^{n}\mu_i l_i)}=
\frac{\Gamma(\gamma+\mu_1 n_0)}{\Gamma(\gamma)}>0.$$

Hence, by virtue of (\ref{eq1}) -- (\ref{eq3}), we get
$$\left |E_{(\mu_1,\mu_2,...,\mu_n),\gamma}(z_1,z_2,...,z_n)\right | $$
$$=
\sum\limits_{k=0}^{\infty}
\sum\limits_{l_1+...+l_n=k}(k;l_1,...,l_n)\frac{|z_1|^{l_1}\cdot...\cdot |z_n|^{l_1}}
{\Gamma(\gamma+\sum_{i=1}^{n}\mu_i l_i)}$$ 
$$\leq C\sum\limits_{k=0}^{\infty}\frac{1}{\Gamma(\gamma+\mu_1 k)}
\sum\limits_{l_1+...+l_n=k}(k;l_1,...,l_n)\prod\limits_{i=1}^{n}|z_i|^{l_i}$$
$$=C\sum\limits_{k=0}^{\infty}\frac{(|z_1|+...+|z_n|)^k}{\Gamma(\gamma+\mu_1 k)}
=C E_{\mu_1,\gamma}(|z_1|+|z_2|+...+|z_n|),$$  
where $C=1+C_0.$ 

The Lemma 1 is proved.


\begin{thebibliography}{99}

\bibitem{gorluch} Y. Luchko, R. Gorenflo, An operational method for solving fractional differential equations with the Caputo derivatives, Acta Math. Vietnam 24 (1999) 207-233.

\bibitem{dzh} M.M. Dzherbashyan, {\it Integral Transforms and Representations of Functions in the Complex Plane}. Nauka, Moscow, 1966 (In Russian).



\end{thebibliography}
\end{document}